\documentclass[12pt]{amsart}
\usepackage{amsfonts,amsthm,amscd}

\newtheorem{theorem}{Theorem}[section]
\newtheorem{corollary}[theorem]{Corollary}
\newtheorem{lemma}[theorem]{Lemma}
\newtheorem{proposition}[theorem]{Proposition}
\theoremstyle{definition}
\newtheorem{defn}{Definition}

\theoremstyle{definition}
\newtheorem{remark}{Remark}
\newtheorem*{acknow}{Acknowledgments}

\newcommand{\nc}{\newcommand}
\newcommand{\be}{\begin{equation}}
\newcommand{\ee}{\end{equation}}
\newcommand{\bea}{\begin{eqnarray}}
\newcommand{\eea}{\end{eqnarray}}
\newcommand{\no}{\nonumber}
\nc{\bth}{\begin{theorem}}
\nc{\eth}{\end{theorem}}
\nc{\bpr}{\begin{proposition}}
\nc{\epr}{\end{proposition}}
\nc{\ble}{\begin{lemma}}
\nc{\ele}{\end{lemma}}
\nc{\bco}{\begin{corollary}}
\nc{\eco}{\end{corollary}}
\nc{\bre}{\begin{remark}}
\nc{\ere}{\end{remark}}

\nc{\f}{\frac}
\nc{\rw}{\rightarrow}
\nc{\Rw}{\Rightarrow}

\begin{document}
\title[Totally real minimal tori in $\mathbb CP^2$]{Totally real minimal tori in $\mathbb CP^2$}
\author{Hui Ma}
\author{Yujie Ma}
\address{                            
    Department of Mathematical Sciences,
    Tsinghua University,
    Beijing 100084, 
    China}
\email{hma@math.tsinghua.edu.cn}
\address{Institute of Systems Science,
Chinese Academy of Sciences,
Beijing 100081,
China}
\email{yjma@mmrc.iss.ac.cn}
\date{}

\begin{abstract}
     In this paper we show that all totally real superconformal minimal tori in $\mathbb CP^2$  correspond with  doubly-periodic finite gap solutions of the Tzitz\'eica equation $$\omega_{z\bar z}=e^{-2\omega}-e^{\omega}. $$
Using the results on the Tzitz\'eica equation in integrable system theory, we describe explicitly all these tori
by Prym-theta functions.
\end{abstract}
\maketitle

\section*{Introduction}

Over the past few years the integrable system approach played an important role in the theory of minimal surfaces 
and harmonic maps. One of the most classical and striking results is on the tori with constant mean curvature (CMC) in    
${\mathbb R}^3$, which correspond the doubly-periodic solutions to the sinh-Gordon equation.
Firstly let us recall the history of the research on CMC tori since it contains all the seeds useful to our case. In 1987, 
Wente's significant existence theorem \cite{Wente1} provoked Abresch \cite{Abresch} to classify all CMC tori having one family of planar curvature lines. These surfaces are given explicitly in terms of  elliptic functions by reducing the PDE to ODEs solvable.  The approach to seek the general solution to the problem has developed into a major area of 
research. The main points of view are due to Pinkall and Sterling \cite{PS} et al, whose approach can be expressed in terms of Hamiltonian systems and Loop groups, Bobenko \cite{Bobenko91},  who sees it in terms of the much studied finite-gap 
solutions for the sinh-Gordon equation and the corresponding Baker-Akhiezer function, and Hitchin \cite{Hit1} whose approach is influenced by twistor theory. 
Later on, Bolton-Pedit-Woodward (\cite{BPW}) showed that the existence of the correspondence between minimal surfaces in $\mathbb CP^n$ and the solutions of affine Toda equations for $SU(n+1)$. And they further proved that any superconformal harmonic $2$-torus in $\mathbb CP^n$ is of finite type. They posed that together with the results in \cite{EW}
this accounts for all (conformal) harmonic $2$-tori in $\mathbb CP^2$. 
Unfortunately, even in the case of $\mathbb CP^2$, the Clifford torus
$$x(z)=[e^{z-{\bar z}}, e^{\epsilon z- {\bar \epsilon}{\bar z}},e^{\epsilon^{2}z-{\bar \epsilon}^{2}{\bar z}}],  \quad \epsilon=e^{\f{2\pi {\bf i}}{3}},$$
which is an isometric embedding of flat torus whose corresponding lattice in $\mathbb C$ is generated by $2\pi/3$ and $2\pi {\bf i}$, is the only example of superconformal harmonic tori in $\mathbb CP^2$ known to us during a long period. On the other hand, there has been much attention given to  totally real (Lagrangian) submanifolds in $\mathbb CP^n$ in the classical differential geometry community. Then a natural problem was posed to give an explicit construction of totally real minimal tori in $\mathbb CP^2$. Moreover, the explicit construction of Lagrangian minimal tori in $\mathbb CP^2$ arouses much interest in the scope of the construction of special Lagrangian 3-folds in $\mathbb C^3$ \cite{Joyce}.
 
In fact, the research on totally real minimal tori in $\mathbb CP^2$ is similar to 
one on CMC tori in $\mathbb R^3$. 
All totally real superconformal minimal tori in $\mathbb CP^2$  correspond with  doubly-periodic finite-gap solutions of
 the  Tzitz\'eica equation 
\begin{equation}\label{db}
\omega_{z\bar z}=e^{-2\omega}-e^{\omega},
\end{equation}
 which is the nearest relative of the sinh-Gordon equation. 
The first progress was given by Castro-Urbano (\cite{CU}). 
They reduced the PDE to ODE solvable by elliptic functions and gave a new family of minimal Lagrangian tori 
in $\mathbb CP^2$, which is characterized by its invariability by a one-parameter group of holomorphic isometrics 
of $\mathbb CP^2$.

On the other hand, the  Tzitz\'eica equation is a soliton equation of physical interest which can be obtained 
as the Toda lattice system of the irreducible nonreduced root system $\it bc_1$. 
Its another form $\omega_{xt}=e^{\omega}-e^{-2\omega}$ was originally derived by Tzitz\'eica in 1910 
as the Gauss-Codazzi equation in asymptotic coordinates $u,v$ for surfaces with the {\it Tzitz\'eica property}. 
It was rediscovered in a solitonic context by Bullough and Dodd in 1977 \cite{DB}. 
Mikhailov proved that this equation is a complete integrable equation in 1980 \cite{Mikhailov}. 
It was also investigated in detail from the point view of integrable systems and affine geometry
by Shabat-Zhiber, Bobenko, etc. It should be mentioned that the finite-gap solutions of this equation was studied by
 Cherdantsev-Sharipov in 1989 \cite{CS}. And the latter studied the CR immersion of minimal tori into
 five-dimensional sphere in $\mathbb C^3$ using the finite-gap solution of the Tzitz\'eica equation \cite{Sharipov}.

In this paper we shall follow the recipe to give all of totally real superconformal minimal tori in $\mathbb CP^2$ by Prym-theta functions using the theory of finite-gap solutions of the  Tzitz\'eica equation. The paper is organized as follows: 
in Section 1 we shall be concerned with the general surface theory in $\mathbb CP^2$, 
and then give the structure equations of totally real minimal surface in $\mathbb CP^2$. 
In Section 2, we reformulate the existence condition of  totally real superconformal minimal surfaces in $\mathbb CP^2$ 
in terms of certain loop groups. In Section 3, we use the method in Section 1 in global studies of totally real 
superconformal minimal tori. In \S4 we show how to construct a class of totally real superconformal minimal immersion 
from $\mathbb R^2$ into $\mathbb CP^2$ by integrating certain Hamiltonian ODE's on finite dimensional subspaces 
of the twisted loop algebra. Following directly the standard scheme, 
we prove that any totally real superconformal minimal torus can be constructed by this method, which is called of finite type. 
In Section 5 and Section 6, we use algebro-geometric methods of nonlinear integrable equations from soliton theory 
to provide an explicit description of these tori in terms of Prym-theta functions.

\section{Theory of surfaces in $\mathbb CP^2$}
Let $\mathbb CP^2$ be the 2-dimensional complex projective space endowed 
with the Fubini-Study metric $h_0$ of constant holomorphic sectional curvature 4. 
Then $h={\mathrm Re} h_0= {\f{1}{2}}(h_0 + {\bar h_0})$ defines a Riemannian metric on $\mathbb CP^2$.
Let $S^5$ be the unit hypersphere in ${\mathbb C^3}$ and let $\pi: S^5 \rightarrow {\mathbb CP^2}$ be the Hopf projection. For any local section $Z:{\mathbb B^4} \rightarrow S^5$ of $\pi$ defined on an open ball
 $\mathbb B^4$ of $\mathbb CP^2$, we have
\begin{equation}\label{fsm}
h_0=(dZ-(dZ \cdot \overline{Z})Z) \otimes
  (d\overline{Z}-(d\overline{Z} \cdot Z)\overline{Z}),
\end{equation}
which we call a local representation of the Fubini-Study metric. 
Here we use the canonical scalar product
$$Z\cdot W=\sum_{k=1}^{3}z_{k}w_{k}.$$
>From \eqref{fsm} we know that the Riemannian metric $h$
on ${\mathbb B^4} \subset {\mathbb CP^2}$ is given by 
\begin{eqnarray*}
h&=&{\f{1}{2}} (h_0 + \overline{h_0})\no\\
  &=&{\f{1}{2}} \{ (dZ-(dZ \cdot \overline{Z})Z)\otimes
   (d \overline{Z} -(d \overline{Z} \cdot Z)\overline{Z})\\
  & &+(d\overline{Z}-(d\overline{Z} \cdot Z)\overline{Z})
   \otimes(dZ-(dZ \cdot \overline{Z})Z)\}.\no
\end{eqnarray*}

Let $x:M\rw {\mathbb CP^2}$ be an immersion of an oriented surface.
The induced metric on $M$ generates a complex structure with respect to which the metric is \begin{equation*}
g=2e^{\omega}dzd{\bar z},
\end{equation*}
where $z$ is a local complex coordinate on $M$ and $\omega : U \rightarrow {\mathbb R}$ is a real function on  an open set $U$ of $M$.
For any local section $Z$ of $\pi : S^5 \rightarrow {\mathbb CP^2}$ we can define a local lift
$y:=Z\circ x$ of the immersion $x: M \rightarrow {\mathbb CP^2}$. Such a local lift $y$ of $x$ exists around each point of $M$.
Let $y: U\rightarrow S^5$ be a local lift of $x$ defined on an open set $U$ of $M$. We denote the Cauchy-Riemannian operators by
$$
 {\frac{\partial}{\partial z}}=
		{\f{1}{2}}({\frac{\partial}{\partial u}}-{\bf i}
		{\frac{\partial} {\partial v}}),\quad\quad\quad
		{\frac{\partial}{\partial {\overline z}}}=
		{\f{1}{2}}({\frac{\partial}{\partial u}}+
	       {\bf i}{\frac{\partial}{\partial v}})
$$
and define
\begin{equation}\label{xi}
 \xi:= y_{z}- ( y_{z} \cdot \overline{y})  y,  \quad\quad 
    \eta:=y_{\overline z}-
        (  y_{\overline z}\cdot {\overline y})  y .  
\end{equation}
The metric $g$ is conformal that gives
\begin{equation}\label{o}
\xi \cdot {\overline \eta} = \xi \cdot {\overline y}  = \eta \cdot {\overline y} =0,  
\end{equation}
\begin{equation}\label{o2}
e^{-\omega}\xi\cdot {\bar \xi}+e^{-\omega}\eta\cdot {\bar \eta}:=a+b=2,	
\end{equation}
where we define
\begin{equation*}
 a:=e^{-\omega}\xi\cdot{\bar\xi},\quad\quad
        b:=e^{-\omega}\eta\cdot{\bar\eta}. 
\end{equation*}
It is easy to check that $a$ and $b$ are independent of the choice of the local lift $y$ and the complex coordinate $z$.
Since $0\leq a,b\leq 2$, we can define globally an invariant  $\theta: M\to [0,\pi]$ as follows
\begin{equation}\label{theta}
 \theta:=2\,{\rm arccos}(\sqrt{\f{a}{2}}). 
\end{equation}
It is easy to verify that the invariant $\theta$ defined above is exactly the K\"ahler angle of $x$. 

\begin{defn}
   A point $p\in M$ is called holomorphic (resp. anti-holomorphic, real) for $x: M\rightarrow {\mathbb C}P^2$ 
if $\theta(p)=0$ (resp. $\pi, {\f{\pi}{2}}$). 
A point is called a complex point of $x$ if it is holomorphic or anti-holomorphic.
\end{defn}

In order to get the complete $SU(3)$-invariant system of $x$ we introduce another two global 
invariants:
\begin{equation}\label{dphi}
\Phi := e^{-\omega}\xi_{\bar z}\cdot {\bar \eta}dz:=\phi dz;
\end{equation}
\begin{equation}\label{dpsi}
\Psi := \xi_{z}\cdot {\bar \eta}dz^{3} :=\psi dz^{3}. 
\end{equation}
Using \eqref{o} we can easily verify that $\Phi$ and $\Psi$ are independent of the choice of local lift $y$ of $x$ and the complex coordinate $z$ of $M$ and thus globally defined on $M$. Moreover, if $x$ is transformed by an isometry $T\in SU(3)$, then $\xi$ and $\eta$ are also transformed by the $T$. It follows by definition that $\Phi$ and $\Psi$ are $SU(3)$ invariant. 
We call $\Psi$ the {\it {cubic Hopf differential}}.
\bre\label{compt}
   Definitions of $\Phi$ and $\Psi$ show that the complex points of $x$ are zeros of $\Phi$ and $\Psi$. 
\ere

Let $x:M\rightarrow {\mathbb CP^2}$ be a surface without complex points and $y: U\rightarrow S^5$ be a local lift of $x$. 
We define $\xi$ and $\eta$ by \eqref{xi}. Then at each point of $U$, $\{y,\xi , \eta \}$ define a basis of $\mathbb C^3$, which, due
to \eqref{o},\eqref{o2} and the fact that $y\cdot {\bar y}=1$, satisfies the following equations:
\begin{equation*}
\sigma_{z}=\sigma {\mathcal U}, \quad \sigma_{\bar z}=\sigma {\mathcal V}, \quad \sigma =(y,\xi, \eta),
\end{equation*}
\begin{eqnarray*}
{\mathcal U}=\left(\begin{array}{ccc}
			\rho& 0& -be^{\omega}\\
			1& (\log ae^{\omega})_{z}+ \rho +a^{-1}\phi & -a^{-1}{\bar \phi}\\
			0& b^{-1}e^{-\omega}\psi &	\rho + b^{-1}\phi
			\end{array}
			\right),
\end{eqnarray*}
\begin{eqnarray*}
{\mathcal V}=\left(\begin{array}{ccc}
	- \rho & -ae^{\omega} & 0\\
	0 & - {\bar \rho}-a^{-1}{\bar \phi} & -a^{-1}e^{-\omega}{\bar \psi}\\
	1 & b^{-1} \phi & (\log be^{\omega})_{\bar z}- {\bar \rho}-b^{-1}{\bar \phi}
	\end{array}\right),
\end{eqnarray*}		
where
\begin{equation}\label{rho}
\rho = y_{z} \cdot {\bar y}.
\end{equation}
The compatibility condition 
$${\mathcal U}_{\bar z}-{\mathcal V}_z =[{\mathcal U},{\mathcal V}]$$
have the following form:
\begin{equation}\label{ie0}
	 \rho_{\bar z}+{\bar
\rho}_{z}=(a-b)e^{\omega},
\end{equation} 
	
\begin{equation}\label{ie1}
(\log(ae^{\omega}))_{z{\bar z}}
	 	=(b-2a)e^{\omega}-(a^{-1}\phi)_{\bar z}-
		(a^{-1}{\bar \phi})_{z}-(ab)^{-1}|\phi|^{2}
		+(ab)^{-1}e^{-2\omega}|\psi|^{2},
\end{equation}	 
\begin{equation}\label{ie2}
\psi_{\bar z}+(a^{-1}-b^{-1}){\bar \phi}\psi=
e^{\omega}(\phi_{z}-\omega_{z}\phi)+(b^{-1}-a^{-1})e^{\omega}\phi^{2} 
-e^{\omega}\phi(\log (ab))_{z},
\end{equation}
\begin{equation}\label{ie3}	
 (\log(be^{\omega}))_{z{\bar z}}=
		(a-2b)e^{\omega}+(b^{-1}\phi)_{\bar
z}+(b^{-1}{\bar \phi})_{z}-(ab)^{-1}|\phi|^{2}+
(ab)^{-1}e^{-2\omega}|\psi|^{2}.
\end{equation}

We note that Formula \eqref{ie0} is not essential. In fact, 
$$(\rho-\rho_0)_{\bar z}+\overline{(\rho-\rho_0)_{\bar z}}=0, \quad
   \rho_0={\f{1}{2}}\int_{M}(a-b)e^{\omega}d{\bar z},$$
which implies that $i\{(\rho-\rho_0)dz- \overline{(\rho-\rho_0)}d{\bar z}\}$ is a closed real 1-form. 
Thus we can write locally this 1-form as $d\vartheta$ for some real function $\vartheta$. 
Since $\rho_0=\rho+i\vartheta_z$, it is easy to check that if we change the local lift $y:U\rw S^5$ to ${\tilde y}=e^{i\vartheta}y$, then the function $\rho$ defined by \eqref{rho} will change to $\rho_0$. Therefore $\{g,\theta,\Phi,\Psi\}$ form a complete $SU(3)$-invariant system of an isometric immersion $x:M\rw {\mathbb CP^2}$ from a simply connected surface $M$ without complex points, i.e. we have the following Fundamental Theorem in \cite{LW}:
\bth[Fundamental Theorem]
\begin{itemize} 
\item[i)] Let $g$ be a Riemannian metric on a simply connected surface $M$, which induces a complex structure on $M$.
Let $\theta:M\rightarrow (0,\pi)$ be a smooth real function. Let $\Phi$ and $\Psi$ be a $(1,0)$-form and a $(3,0)$-form on 
$M$ respectively. If $\{g,\theta,\Phi,\Psi\}$ satisfy the relation \eqref{ie1}, \eqref{ie2} and \eqref{ie3}, then there exist an 
immersion $x:M\rightarrow {\mathbb CP^2}$ without complex point such that $g$ is its induced metric, and 
$\{\theta, \Phi,\Psi\}$ are exactly the $SU(3)$-invariants of $x$ defined by \eqref{theta},\eqref{dphi} and \eqref{dpsi}.
\item[ii)] If two surfaces $x,{\tilde x}:M\rw {\mathbb C}P^2$ without complex points have the same invariant system, 
i.e. there exists a diffeomorphism $\sigma : M\rw M$ such that $g=\sigma^{*} \tilde{g}, \theta= \tilde{\theta}\circ \sigma,  \Phi = \sigma^{*} {\tilde \Phi},
\Psi=\sigma^{*} \Psi$, then there is an isometry $T\in SU(3)$ with  $x=T \circ {\tilde x} \circ \sigma$.
\end{itemize}
\eth
\bre\label{min}   
   It is not hard to show that $x:M\rw {\mathbb CP^2}$ is minimal if and only if $\Phi\equiv 0$.
\ere

Suppose that $x:M\rw {\mathbb CP^2}$ is minimal, then the Gauss-Codazzi equations become
\begin{eqnarray*}
(\log (ae^{\omega}))_{z{\overline z}}&=&(ab)^{-1}e^{-2\omega}|\psi|^2+(b-2a)e^{\omega},\\
                               \psi_{\overline z} &=&0,\\
 (\log (be^{\omega}))_{z{\overline z}}&=&(ab)^{-1}e^{-2\omega}|\psi|^2+(a-2b)e^{\omega},
   \end{eqnarray*}
which are invariant with respect to the transformation
$$\Psi\rw {\lambda^*}\Psi,\quad |\lambda^*|=1.$$

\bpr
If $x:M\rw {\mathbb CP^2}$ is a minimal surface, then the Hopf differential $\Psi$ is a holomorphic cubic form on $M$. In particular, if $\Psi \not\equiv 0$, then the complex points of $x$ are isolated.
\epr
\begin{defn}
A minimal surface $x:M\rw {\mathbb CP^2}$ without complex points is called {\it superminimal}
if $\Psi \equiv 0$, otherwise, {\it superconformal} if $\Psi \not\equiv 0$.  
\end{defn}
 
Treating $\lambda^*$ as a deformation parameter we obtain 
\bth\label{family}
Every superconformal minimal surface  has a one-parameter family of deformations preserving the induced metric and the K\"ahler angle. 
\eth

>From now on we suppose that $x: M\rw {\mathbb CP^2}$ is a totally real minimal surface. Thus we further have  $\Phi \equiv 0$ and $a=b\equiv 1$. 
Without loss of generality a change of local lift confines ourselves to the case
 $\rho \equiv 0$. Then we 
get 
\begin{equation}\label{e1}
       \sigma_{z}=\sigma {\mathcal U}, \quad \sigma_{\bar z}=\sigma {\mathcal V}, \quad \sigma =(y,\xi, \eta),
\end{equation}
\begin{eqnarray}\label{e2}
{\mathcal U}=\left(\begin{array}{ccc}
			0& 0& -e^{\omega}\\
			1& \omega_{z} &0\\
			0&e^{-\omega}\psi &0
			\end{array}
			\right),
			{\mathcal V}=\left(\begin{array}{ccc}
				0 & -e^{\omega} & 0\\
				0 & 0 & -e^{-\omega}{\bar \psi}\\
				1 & 0 & \omega_{\bar z}
				\end{array}\right).
\end{eqnarray}
In this case the Gauss-Codazzi equations become	
\begin{equation}\label{e3}
\omega_{z{\bar z}}=e^{-2\omega}|\psi|^2-e^{\omega}, \quad \psi_{\bar z}=0.
\end{equation}

The framing $\sigma(z,{\bar z},\lambda^*)$ solving the system  \eqref{e1} with
\begin{eqnarray}
{\mathcal U}(\lambda^*)=\left(\begin{array}{ccc}
			0& 0& -e^{\omega}\\
			1& \omega_{z} &0\\
			0&{\lambda^*} e^{-\omega}\psi &0
			\end{array}
			\right),
			{\mathcal V}(\lambda^*)=\left(\begin{array}{ccc}
				0 & -e^{\omega} & 0\\
				0 & 0 & -{\bar {\lambda^*}}e^{-\omega}{\bar \psi}\\
				1 & 0 & \omega_{\bar z}
				\end{array}\right)\no
\end{eqnarray}	
describes a family of totally real superconformal minimal surfaces $x^{\lambda^*}:M\rw {\mathbb CP^2}$.

In order to simplify the framing, we pass to a gauge equivalent frame function
\begin{eqnarray}
F(\nu)=\sigma(\lambda^*)\left(\begin{array}{ccc}
                                                 1&0&0\\
                                                 0&{\f{1}{{\bf i}\nu}}e^{-\f{\omega}{2}}&0\\
                                                 0&0&-{{\bf i}\nu}e^{-\f{\omega}{2}}
                                                   \end{array}
                                               \right)\in SU(3)\no
\end{eqnarray}
where 
\begin{equation}\label{norm}
{\nu}^{3}={\bf i}{\lambda^*}.
\end{equation}

Combined with Theorem \ref{family}, we obtain
\bth\label{t:3.4} 
Let $F(z,{\bar z},\nu), \nu\in S^1$ be a solution of the system 
\begin{equation}\label{stan}
F_z = F U(\nu), \quad F_{\bar z} = FV(\nu),
\end{equation}
\begin{eqnarray}\label{1}
U(\nu)=\left(\begin{array}{ccc}
0&0&{\bf i}\nu e^{\f{\omega}{2}}\\
{\bf i}\nu e^{\f{\omega}{2}}&\f{\omega_{z}}{2}&0\\
0&-{\bf i}\nu \psi e^{-\omega}&-{\f{\omega_{z}}{2}}
\end{array}
\right),
\end{eqnarray}
\begin{eqnarray}\label{2}
V(\nu)=\left(\begin{array}{ccc}
0&{\f{\bf i}{\nu}}e^{\f{\omega}{2}}&0\\
0&-{\f{\omega_{\bar z}}{2}}&-{\f{\bf i}{\nu}}{\bar \psi}e^{-\omega}\\
{\f{\bf i}{\nu}}e^{\f{\omega}{2}}&0&{\f{\omega_{\bar z}}{2}}
\end{array}
\right),
\end{eqnarray}
normalized by \eqref{norm}
where $\psi \neq 0$, thus $[F(\nu)(e_0)]$ gives a totally real superconformal minimal surface in 
$\mathbb CP^2$ with metric $g=2e^{\omega}dzd{\bar z}$, and nonzero Hopf differential 
$\Psi^{\lambda^*}={\lambda^*}\Psi={\lambda^*}\psi dz^3$.

Conversely, suppose $x^{\lambda^*}:M\rw {\mathbb CP^2}$ is a conformal parameterization of a totally real 
superconformal minimal surface in $\mathbb CP^2$ with the metric $g=2e^{\omega}dzd{\bar z}$ and nonzero Hopf differential 
$\Psi^{\lambda^*}={\lambda^*} \psi dz^3$. Then there exists a unique frame $F^{\lambda^*}: M\rw SU(3)$
is a solution of \eqref{1},\eqref{2} as above.    
\eth

\bre
   In a neighborhood of a point $\psi\neq 0$ by a conformal change of coordinate $z\rw w(z)$ one can always normalize $\psi=-1$. Thus the Gauss equation becomes the Tzitz\'eica equation 
$\omega_{z\bar z}=e^{-2\omega}-e^{\omega}$.
\ere

\section{Loop group formulation}

The above results may be conveniently described in terms of certain loop groups.
Put $G=SU(3)$ and $G^{\mathbb C}=SL(3,{\mathbb C})$ denotes the complexification of $G$.
The corresponding Lie algebras are denoted by $\mathcal G = su(3)$ and ${\mathcal G}^{\mathbb C}=su(3,{\mathbb C})$, respectively.

Let $\epsilon=e^{2\pi {\bf i}/3}$. Put
$$Q=\left(\begin{array}{ccc}
			1& 0& 0\\
			0 & \epsilon & 0\\
			0 & 0 & \epsilon^2
			\end{array}
			\right), \quad
 T=\left(\begin{array}{ccc}
		1 & 0 & 0\\
		0 & 0 & 1\\
		0 & 1 & 0
	\end{array}\right).
 			$$		
An inner automorphism  $\gamma : G^{\mathbb C}\rw G^{\mathbb C}$ of order $3$ is defined by
$g \mapsto  Q g Q^{-1}$
and an involutive automorphism 
$\sigma : G^{\mathbb C}\rw G^{\mathbb C}$ is defined by $g  \mapsto T(g^{t})^{-1}T$.
The corresponding automorphism  $\gamma$ of order $3$  
and involutive automorphism $\sigma$ of ${\mathcal G}^{\mathbb C}$ are 
$$\gamma : {\mathcal G}^{\mathbb C} \to {\mathcal G}^{\mathbb C}, \quad  \xi \mapsto  Q\xi Q^{-1},$$
$$\sigma : {\mathcal G}^{\mathbb C}\to {\mathcal G}^{\mathbb C}, \quad \xi  \mapsto -T\xi^{t}T.$$

We denote by ${\mathcal M}_k$ the  $\epsilon^k$-eigenspace of $\gamma$. Explicitly, they are described as
$${\mathcal M}_{0}=\left\{
                               \left(\begin{array}{ccc}
			* & 0& 0\\
			0 & * & 0\\
			0 & 0 & *
			\end{array}
\right)\in {\mathcal G}^{\mathbb C}
 \right \},$$
$${\mathcal M}_{1}=\left\{
                   \left(\begin{array}{rcc}
			0 & 0& *\\
			{*} & 0 & 0\\
			0 & * & 0
			\end{array}
\right)\in {\mathcal G}^{\mathbb C}
 \right \},  $$
$${\mathcal M}_{2}=\left \{
            \left(\begin{array}{rcc}
			0 & * & 0\\
			0 & 0 & *\\
			{*} & 0 & 0
			\end{array}
	\right)\in {\mathcal G}^{\mathbb C}
           \right \}.$$
Obviously, it holds that
$$[{\mathcal M}_{k},{\mathcal M}_{j}]\subset {\mathcal M}_{k+j},\quad {\mathcal M}_{k+3j}={\mathcal M}_{k}.$$

We denote by the ${\mathcal N}_k$ $(-1)^k$-eigenspace of $\sigma$. Explicitly, they are described as
$${\mathcal N}_{0}=\left\{
                         \left(\begin{array}{ccc}
			0& b_1& b_2\\
			-b_2 & a & 0\\
			-b_1 & 0 & -a
			\end{array}
			\right)\in {\mathcal G}^{\mathbb C};a,b_i \in {\mathbb C},i=1,2
                               \right \},  $$
$${\mathcal N}_{1}=\left\{
                     \left(\begin{array}{ccc}
			2a& b_1& b_2\\
			b_2 & -a & b_3\\
			b_1 & b_4 & -a
			\end{array}
	\right)\in {\mathcal G}^{{\mathbb C}};a,b_i \in {\mathbb C},i=1,\cdots ,4
                                   \right \}. $$
Similarly, we have
$$[{\mathcal N}_{k},{\mathcal N}_{j}]\subset {\mathcal N}_{k+j}, \quad {\mathcal N}_{k+2j}={\mathcal N}_{k}.$$

The matrices
$$A=U(\nu)+V(\nu), \quad B={\bf i}(U(\nu)-V(\nu))$$
corresponding to real vector fields $\partial_{x}=\partial_{z}+\partial_{\bar z}$ and 
$\partial_{y}={\bf i}(\partial_{z}-\partial_{\bar z})$ belong to the twisted loop algebra
$$\Lambda {\mathcal G}^{\mathbb C}_{\gamma,\sigma}=\{\xi : S^{1}\rw  {\mathcal G}^{\mathbb C}\;|\;        \gamma(\xi(\nu))=\xi(\epsilon \nu),
          \sigma(\xi(\nu))=\xi(- \nu)
               \},$$
and $F$ in \eqref{stan} lies in the corresponding twisted loop group
$$\Lambda G^{\mathbb C}_{\gamma,\sigma}=\{ g: S^{1}\rw G^{{\mathbb C}}\;|\;
    \gamma(g(\nu))=g(\epsilon \nu),
\sigma(g(\nu))=g(- \nu)
   \}.$$
Conjugation on $\Lambda {\mathcal G}^{\mathbb C}_{\gamma,\sigma}$ is defined as 
$${\bar \xi}(\nu):={\overline {\xi(1/{\bar \nu})}}.$$
We put
$$\Lambda {\mathcal G}_{\gamma,\sigma}=\{ \xi\in \Lambda {\mathcal G}^{\mathbb C}_{\gamma,\sigma}|
          \bar{\xi}(\nu)=\xi(\nu)
                        \}$$
and denote the corresponding loop group by $\Lambda G_{\gamma,\sigma}$.
$\Lambda G^{\mathbb C}_{\gamma,\sigma}$ is a Banach Lie group and $\Lambda {\mathcal G}^{\mathbb C}_{\gamma,\sigma}$ admits the ${\hbox{Ad}}\Lambda G^{\mathbb C}_{\gamma,\sigma}$-invariant inner product
$$(\xi,\eta)_{L^2}=\int_{S^1}(\xi(\nu),\eta(\nu))_{{\mathcal G}^{\mathbb C}}d\nu,$$
where $(\cdot,\cdot)_{{\mathcal G}^{\mathbb C}}$ is the ad-invariant inner product 
on ${\mathcal G}^{\mathbb C}$.
Expanding $\xi \in \Lambda{\mathcal G}^{{\mathbb C}}_{\gamma,\sigma}$ in a Laurent series we get 
$$\xi = \sum_{k\in {\mathbb Z}}\nu^{k}\xi_{k}$$ 
with $\xi_{k}\in {\mathcal M}_{k} \cap {\mathcal N}_{k}$. 
Further for any $\xi \in \Lambda{\mathcal G}_{\gamma,\sigma}$ we have $\xi_{-k}={\bar \xi}_k$. 

We introduce a filtration of $\Lambda{\mathcal G}_{\gamma,\sigma}$ by finite-dimensional subspaces
$$\Lambda_1 \subset \Lambda_2\subset \cdots \subset \Lambda {\mathcal G}_{\gamma,\sigma},$$
 setting
$$\Lambda_{d}=\{\xi\in \Lambda{\mathcal G}_{\gamma,\sigma}|\xi = 
\sum_{|k| \le d}\nu^{k}\xi_{k}
                                          \}$$
for $d\in {\mathbb N}$.
\begin{defn}
$\theta (\nu) = (\nu^{-1}V_{1}+V_{0})d{\bar z}+ (U_0+\nu U_1)dz$ is 
called a family of {\it real normalized admissible connections} if it satisfies 
\begin{itemize}
\item[(i)] $\theta(\nu)\in \Lambda_1$,
\item[(ii)] The connections $d+\theta$ are flat for arbitrary $\nu\in S^{1}$,
\item[(iii)] $\hbox{Tr} (U_1)^3=-3{\bf i}$.
\end{itemize}
\end{defn}

For any totally real superconformal minimal surface $x: M\rw {\mathbb CP^2}$ with $\psi=-1$,  there exists a family of frames $F^{\nu}$ such that the connections $\theta(\nu)=(F^{\nu})^{-1}dF^{\nu}$
having the form as \eqref{1}, \eqref{2} with $\psi=-1$ are real normalized admissible.

Conversely, given a family of real normalized admissible connection $\theta\in \Lambda_1$,
there exists a family of frames $F: M\rw \Lambda G_{\gamma,\sigma}$ such that $F^{-1}dF=\theta$, that is, 
$$F^{-1}F_{z} = U_0 + \nu U_1=
               \left(\begin{array}{ccc}
                       0& & \\
                        &a& \\
                      & &-a
                               \end{array}
                                    \right)
                         + \nu \left(\begin{array}{ccc}
                                   0&0&b_0\\
                                   b_0&0&0\\
                                   0&b_1&0
                           \end{array}
                   \right),
$$
and $F^{-1}F_{\bar z}={\overline {(F^{-1}F_z)}}$, where $a,b_0,b_1\in {\mathbb C}$.
Since 
$$U_{1}=\left(\begin{array}{ccc}
0&0&b_0\\
b_0&0&0\\
0&b_1&0
\end{array}
\right) \in {\mathcal M}_{1} \cap {\mathcal N}_{1}, \quad {\hbox {Tr}}(U_{1})^{3}=-3{\bf i},$$
and ${\hbox {Tr}}(U_{1})^{3} = 3b_{0}^{2} b_{1}$, we get $b_{0}^{2} b_{1} =-{\bf i}$. Assume that
$$b_0 = {\bf i}e^{\f{\omega+ {\bf i} \theta}{2}}, 
\quad  b_1 ={\bf i}e^{-\omega-{\bf i} \theta}, \quad  \omega,\theta \in {\mathbb R}.$$
Thus there exists 
$$W=\left(\begin{array}{ccc}
0& & \\
 &{\f{\omega + {\bf i} \theta}{2}}& \\
 & &{\f{-\omega- {\bf i}\theta}{2}}
\end{array}
\right) \in {\mathcal M}_{0} \cap {\mathcal N}_{0},$$
such that 
\begin{equation}\label{*}
e^{W}Be^{-W}=U_{1},
\end{equation}
where
$$B=\left(\begin{array}{ccc}
0&0&{\bf i}\\
{\bf i}&0&0\\
0&{\bf i}&0
\end{array}
\right).$$
Since $M$ is simply connected, $W$ is defined on $M$.
By \eqref{*}, $U_{1{\bar z}}= [W_{\bar z}, U_{1}]$. On the other hand, from the facts that 
$U_{1{\bar z}}=[U_{1}, {\overline{U_0}}],$ and ${\hbox {ad}}U_{1}$ is injective on  
${\mathcal M}_{0} \cap {\mathcal N}_{0}$, we have
$$W_{\bar z}=-{\overline{U_0}}.$$
Put
$$W=\left(\begin{array}{ccc}
0& & \\
 &{\f{\omega}{2}}& \\
 & &-{\f{\omega}{2}}
\end{array}
\right)
+ \left(\begin{array}{ccc}
0& & \\
 &{\bf i}{\f{\theta}{2}}& \\
 & &-{\bf i}{\f{\theta}{2}}
\end{array}
\right):=\Omega + \Lambda,$$
thus $-{\bar \Omega}=\Omega, {\bar \Lambda}=\Lambda$.

Using the commutativity of ${\mathcal M}_0$, by a direct computation we have 
$$F^{-1}F_{\bar z} = - \Omega_{\bar z} - \Lambda_{\bar z} 
       + \nu^{-1}e^{\Lambda}e^{-\Omega}{\bar B}e^{\Omega}e^{- \Lambda}.$$
Passing a gauge transformation ${\hat F} = Fe^{\Lambda}$, we obtain
\begin{equation*}
{\hat F}^{-1}{\hat F}_{\bar z}=-\Omega_{\bar z }+ \nu^{-1}e^{-\Omega}{\bar B}e^{\Omega}.
\end{equation*}
>From the reality conditions it follows that 
$$\begin{array}{rcl}
{\hat F}^{-1}{\hat F}_{z}&=&\Omega_{z}+ \nu e^{\Omega}Be^{-\Omega}\\
&=&
 \left( \begin{array}{ccc}
0&0&{\bf i}\nu e^{\f{\omega}{2}}\\
{\bf i}\nu e^{\f{\omega}{2}}&{\f{\omega_{z}}{2}}&0\\
0&{\bf i}\nu e^{-\omega}&-{\f{\omega_{z}}{2}}
\end{array} 
\right).
\end{array}
$$
Thus ${x^{\nu}}=[{\hat F}(e_0)]$ is a family of totally real superconformal minimal surfaces with the induced metric $g=2e^{\omega}dzd{\bar z}$, nonzero Hopf differential $\Psi^{\lambda^*}=-{\lambda^*} dz^{3}$, where $x^1$ is a normalized totally real superconformal minimal surface. Then we have the following description
in terms of loop groups:
\bth\label{1:1}
There exists a natural correspondence between the  totally real superconformal minimal surface $x: M\rw 
{\mathbb CP^2}$  with $\psi=-1$ and the real normalized admissible connection on $M$.
\eth

\section{Totally real superconformal minimal tori}

Now suppose that $M$ is a compact Riemann surface. If $M$ is of genus $0$ then the holomorphic cubic differential $\Psi=\psi dz^3$ must vanish identically.
That is to say that any totally real minimal 2-sphere in $\mathbb CP^2$ must be superminimal. The classification of totally real minimal tori is not so simple as that of minimal spheres but analytic tools enable us to achieve success in this case just like CMC tori in $\mathbb R^3$. Any Riemann surface of genus one is conformally equivalent to the quotient ${\mathbb C}/\Lambda$ of the complex plane by a lattice $\Lambda$. The corresponding conformal parameterization 
on a torus is given by a doubly-periodic mapping $x: {\mathbb C}/ \Lambda \rw {\mathbb CP^2}$. 
The metric $\omega(z,\bar z)$ and the Hopf differential $\Psi(z,\bar z)$ in this parameterization are  doubly-periodic with respect to the lattice $\Lambda$. Note that $\psi_{\bar z}=0$ and $\psi(z)$ is a bounded elliptic function, thus a constant. This 
constant is not zero, otherwise, 
the Gauss equation becomes the Liouville equation
$$\omega_{z{\bar z}} + e^{\omega}=0, $$
and the Gauss curvature of $M$ is $K=-\omega_{z\bar z}e^{-\omega}=1$, which corresponds to the 
totally real immersion of sphere.
 As before we normalize equation \eqref{e3} to the Tzitz\'eica equation \eqref{db} by $\psi=-1$.
    
Denoting the generators of $\Lambda$ by
$$Z_{1}=X_{1}+{\bf i}Y_{1}, \quad Z_{2}=X_{2}+{\bf i}Y_{2},$$
one obtains the following     
  
\bpr
Any totally real superconformal minimal torus can be conformally parameterized by a doubly-periodic immersion
$F: {\mathbb C}\rw SU(3)$ 
$$F(z+Z_{i}, {\bar z}+{\bar Z}_{i}) = F(z, {\bar z}), i=1,2$$
with $\psi = -1$. In this parameterization the metric $\omega(z,{\bar z})$ of $[F(e_0 )]$ is a doubly-periodic 
solution to the  \eqref{db}.
\epr

To describe all totally real superconformal minimal tori one should solve the following problems.
\begin{itemize}
\item[({\it i})] Describe all doubly-periodic solutions $\omega (z,{\bar z})$ of the equation \eqref{db};
\item[({\it ii})] Find $F(z,{\bar z},\nu)$;
\item[({\it iii})] The first column of $F(z,{\bar z}, \nu)$ describes the corresponding totally real 
superconformal minimal
immersion. In general, this immersion is not doubly-periodic. One should specify parameters of the 
solution $\omega(z,{\bar z})$, which yield doubly-periodic $F(z,{\bar z})$.
\end{itemize}

These three problems will be solved in Section 7 simultaneously using methods of the finite-gap integration theory.

\section{Polynomial Killing field}

In Section 2 we express the condition of a totally real superconformal minimal surface in $\mathbb CP^2$ as the flat condition on a loop algebra valued 1-form with an algebraic constraint (i.e. taking value in $\Lambda_1$). Such a reformulation is the starting point for their integration.
In this section firstly we shall show how to construct real normalized admissible 
$\Lambda_1$-valued 1-forms on the covering space $\mathbb R^2$ of $T^2$ based on a general result on $R$-matrices and commuting flows and characterize the connections obtained by the existence of polynomial 
Killing fields.

Define an $R$-matrix on $\Lambda{\mathcal G}^{\mathbb C}_{\gamma,\sigma}$ by
$$\nu^{k}\xi_{k} \mapsto {\f{1}{2}}\hbox {sign}(k)\nu^{k}\xi_{k}.$$
This definition makes $R$ purely imaginary:
$${\overline {R \xi}} = -R{\bar \xi}.$$

To apply the standard recipe for the construction of $\Lambda_1$-valued flat 1-form on $\mathbb R^2$ we first fix $d\in {\mathbb N}$ with $d \equiv 1\bmod 6$. Define 
$f\in {\mathcal C}^{\infty}(\Lambda{\mathcal G}_{\gamma,\sigma})$ by
$$f(\xi) = \int_{S^1} \nu ^{1-d} (\xi ,\xi)_{\mathcal G} = f_{1}(\xi) - {\bf i}f_{2}(\xi).$$
It is easy to see that $f_{1}, f_{2}$ are Ad$G$-invariant and their gradients with respect to the inner product on 
$\Lambda{\mathcal G}_{\gamma,\sigma}$ are
$$\nabla f_{1}(\xi) - {\bf i}\nabla f_{2}(\xi) =2 \nu^{1-d}\xi.$$

Define the Hamiltonian vector fields $X_{1},X_{2}$ corresponding to $f_{1}, f_{2}$ by
$${\f{1}{2}}(X_{1} - {\bf i}X_{2})(\xi) = [\xi, (R+{\f{1}{2}}) \nu^{1-d}\xi].$$
Since $\xi$ is real, 
$${\f{1}{2}}(X_{1} + {\bf i}X_{2})(\xi) = [\xi, -(R-{\f{1}{2}}) \nu^{d-1}\xi].$$
Restricting $\xi\in \Lambda_d$ these become 
$${\f{1}{2}}(X_{1} - {\bf i}X_{2})(\xi) = [\xi, {\frac{\xi_{d-1}}{2}} + \nu\xi_{d}],$$
$${\f{1}{2}}(X_{1} + {\bf i}X_{2})(\xi) = [\xi, {\frac{\xi_{1-d}}{2}} + \nu^{-1}\xi_{-d}],$$
since $R\pm {\f{1}{2}}$ annihilate 
$$\Lambda^{-} = \{ \xi \in \Lambda {\mathcal G}^{{\mathbb C}}: \xi_{n}=0 \quad {\hbox {for}} \quad n < 0 \}$$
 and  
$$\Lambda^{+} = \{ \xi \in \Lambda {\mathcal G}^{{\mathbb C}}: \xi_{n}=0 \quad {\hbox {for}} \quad n > 0 \}$$
respectively. 
Denote $Z={\f{1}{2}} ( X_{1}-{\bf i}X_{2} )$. 
It is clear that $X_1$ and $X_2$ are tangent to $\Lambda_d$. It follows from Theorem 2.1 of \cite{BFPP} that
$X_i$ are complete and commutative 
so that 
the system of ODEs
\begin{eqnarray}\label{a}
{\f{\partial \xi}{\partial z}}&=&Z(\xi)\no\\
 \frac{\partial \xi}{\partial {\bar z}}&=&{\bar Z}(\xi)\\
 \xi(0)&=&{\stackrel{\circ } {\xi}} \in \Lambda_{d}\no
\end{eqnarray}  
has a unique, global, real solution $\xi: {\mathbb R^2}\rw \Lambda_d$ for any given real initial condition
${\stackrel{\circ } {\xi}}\in \Lambda_d$.

Tr$\xi_{d}^{3}$ is a constant on $\mathbb R^2$. In fact, from
$${\frac{\partial \xi_d}{\partial  z}}={\f{1}{2}}[\xi_{d},\xi_{d-1}], \quad
{\frac{\partial \xi_d}{\partial {\bar z}}}={\f{1}{2}}[\xi_{d},\xi_{1-d}],$$
one gets
$$d(\hbox{Tr}\xi_{d}^{3}) = 
3\hbox{Tr} (\xi_{d}^{2}d{\xi_{d}})
 = 0.$$  
Therefore if assume that 
$\hbox{Tr}({\stackrel{\circ}{\xi}_{d}})^{3} \neq 0$, then 
$\hbox{Tr}\xi_{d}^{3}=\hbox{Tr}({\stackrel{\circ}{\xi}_{d}})^{3}\neq 0$. 

Given a solution $\xi :{\mathbb R^2}\rw \Lambda_d$ of \eqref{a} the $\Lambda_1$-valued 1-form $\theta$ on 
$\mathbb R^2$:
\begin{eqnarray}\label{eq}
{\theta}&=&(\f{\xi_{d-1}}{2} + \nu \xi_{d})dz + ({\f{\xi_{1-d}}{2}} 
+ \nu^{-1}\xi_{-d})d{\bar z}\\
&:=& (U_{0}+\nu U_{1})dz+({\overline U}_0 +\nu^{-1}{\overline U}_1)d{\bar z}\no
\end{eqnarray}
satisfies the Maurer-Cartan equation and thus gives rise to a map $F:{\mathbb R^2}\rw \Lambda G_{\gamma,\sigma}$ with $F^{-1}dF=\theta$. It is easy to see that such a family of connections constructed above is the desired real normalized 
admissible connections. 
Following existing nomenclature we shall give 
\begin{defn}
The map $F: {\mathbb R^2}\rw \Lambda G_{\gamma,\sigma}$ constructed above is called of {\it finite type}.
The corresponding immersion $x:{\mathbb R^2}\rw {\mathbb CP^2}$ is also called of {\it finite type}. 
\end{defn}
By \eqref{eq}, the equation \eqref{a} can be written as
\begin{equation}\label{con}
d\xi=[\xi,\theta].
\end{equation}
\begin{defn}
We call a vector field $\xi : {\mathbb R^2}\rw \Lambda_d$ {\it a polynomial Killing field} of a real normalized admissible connection $\theta=Udz+{\bar U}d{\bar z}$, if $\xi$ satisfies \eqref{con}. Moreover, if $\xi$ satisfies 
$$U={\frac{\xi_{d-1}}{2}} + \nu \xi_{d},$$
we call $\xi$ {\it an adapted polynomial Killing field}.
\end{defn}

We now show that for any totally real superconformal minimal torus in $\mathbb CP^2$ there always exists a polynomial Killing field, that is to say
\bth\label{t6.1}
Any totally real superconformal minimal torus in $\mathbb CP^2$ is of finite type.
\eth
\begin{proof}
 From Theorem \ref{t:3.4}, for any immersion $x:T^2\rw {\mathbb CP^2}$, there exists a unique frame $F^{\nu}:{\mathbb R^2}\rw SU(3)$ such that 
$$(F^{\nu})^{-1}{\frac{\partial F^{\nu}}{\partial z}} = \Omega_{z}+ \nu \hbox{Ad}\exp (\Omega)B$$
where $\Omega : T^2 \rw {\mathcal M_0}\cap {\mathcal N}_0$. 
It is sufficient to construct an adapted polynomial
Killing field $\xi: {\mathbb R^2}\rw \Lambda_{d}$, $d\equiv 1\bmod 6$ such that  
\begin{equation}\label{e6.1} 
d\xi = [\xi, (\Omega_{z}+ \nu \hbox{Ad}\exp \Omega(B))dz 
                    + (-\Omega_{\bar z}+ \nu^{-1}\hbox{Ad}\exp (-\Omega)(\bar B))d{\bar z}]
\end{equation}
There exists a regular algebraic description of these polynomial Killing fields through a formal Killing field 
(see \cite{FPPS},\cite{BPW}), which is in our case  a $\Lambda {\mathcal G}^{\mathbb C}_{\gamma,\sigma}$-valued formal power series solution $Y=\sum_{k \le j } \nu^{k}Y_{k} , j \equiv 1 \, \bmod\, 6$
of $dY = [Y,{\hat F}^{-1}d{\hat F}]$,
where
$${\hat F}^{-1}d{\hat F} = (2\Omega_{z} + \nu B)dz+ \nu^{-1}e^{-2\Omega}{\bar B}e^{2\Omega}d{\bar z}.$$

Since for all $n\in{\mathbb N}$,
$$Y^{(n)}(\nu)=\nu^{6n}Y(\nu)$$
is again a formal Killing field in $\Lambda{\mathcal G}^{\mathbb C}_{\gamma,\sigma}$. 
The diagonal terms
$${\f{1}{2}}Y_{-6n} = u_{n}\left( \begin{array}{rcl}
          0& & \\
            &1& \\
         & &-1  
\end{array}\right)\in {\mathcal M}_{0} \cap {\mathcal N}_{0},\quad n=1,\cdots$$
solve the elliptic equation
\begin{equation}\label{e6.11}
(\partial_{z{\bar z}}+e^{\omega}+2e^{-2\omega})u_{n} = 0.
\end{equation}
\begin{lemma}\label{l6.6}
If $\omega$ is a doubly-periodic solution of \eqref{db}:
$$\omega(z+Z_{i}, {\bar z} + {\bar Z}_{i}) = \omega(z,{\bar z}), \, 
i= 1,2, \quad Im Z_{1}/Z_{2}\neq 0,$$
then only the finite vector fields ${\f{1}{2}}Y_{-6n}, n=1,\cdots$ are linearly independent.
\end{lemma}
\begin{proof}
All $u_n$ are also doubly-periodic. The equation \eqref{e6.11} determines an elliptic linear operator $L$ on the torus $T^2$:
\begin{equation*}
Lu_{n}=(\partial_{z\bar z}+e^{\omega}+2e^{-2\omega})u_{n}=0.
\end{equation*}
By the linear elliptic theory, the compactness of the torus implies
that the spectrum of this operator is discrete and all the eigenspaces 
are finite dimensional. All vector fields ${\f{1}{2}}Y_{-6n}=u_{n}\hbox{diag} (0,1,-1)$ belong to the kernel of $L$.
This observation proves the lemma.
\end{proof}
Thus it follows from the proof of Theorem 3.6 in \cite{BPW} that 
such a polynomial Killing field $\xi$ exists.
Now it is sufficient to prove the existence of a formal Killing field 
$Y(\nu)=\sum_{k\leq j}\nu^{k}Y_{k}$, $j\equiv 1\bmod 6$ of $\hat F$.
Following the standard procedure in \cite{FPPS} we get a formal vector field 
${\hat Y}=\sum_{k=-1}^{\infty} \nu^{-k}{\hat Y}_{-k}$
with values in $\Lambda {\mathcal G}_{\gamma}^{\mathbb C}$ satisfying $d{\hat Y}=[{\hat Y},{\hat F}^{-1}d{\hat F}]$.
Since $\sigma({\mathcal M}_{-k})\subset {\mathcal M}_{-k}$,
$${\tilde Y}=\sum_{k=-1}^{\infty}\nu^{-k}(-1)^{k}\sigma({\hat Y}_{-k})$$
also satisfies the same equation.
 Setting
$$Y={\f{1}{2}}({\hat Y}+{\tilde Y}),$$
we obtain a formal Killing field of $\hat F$ as required, which completes the proof of Theorem \ref{t6.1}. 
\end{proof}

\section{Spectral curve and Baker-Akhizer function}

Let $\omega(z,{\bar z})$ be a solution of the Tzitz\'eica equation \eqref{db} with the 
polynomial Killing field $\xi(\nu)$. The curve
\begin{equation}\label{spec}
\det(\xi(\nu)-\mu I)=0
\end{equation} 
is independent of $z$ and ${\bar z}$, which is called the spectral curve of the solution $\omega(z,{\bar z})$. The matrix $\xi$ is traceless so we see that the spectral curve
is of the form 
$$\mu^3 = {\f{1}{2}}\mu \mbox{tr} \xi^{2}(\nu)+\det \xi(\nu).$$
We shall suppose that the spectral curve is a nonsingular curve  below.

Compactified at $\mu=\infty$ the spectral curve \eqref{spec} determines a compact Riemann surface $\stackrel{\circ}{\Gamma}$. Due to symmetries of the loop algebra 
$\Lambda{\mathcal G}_{\gamma,\sigma}$ beside the 3-order automorphism 
$$\gamma:(\mu,\nu)\rw (\mu, \epsilon \nu),$$
it possesses two more involutions: a holomorphic 
\begin{equation}\label{holo}
\sigma : (\mu,\nu) \rw (-\mu,-\nu)
\end{equation}
and an anti-holomorphic 
\begin{equation}\label{anti}
\tau : (\mu,\nu) \rw (-{\bar \mu},\f{1}{\bar \nu}).
\end{equation}
These symmetries reflect on the expression of the spectral curve $\stackrel{\circ}{\Gamma}$ is 
$${\f{1}{2}}\mbox{tr}\xi^{2}(\nu) = \sum_{k=-\f{d-1}{3}}^{\f{d-1}{3}}p_{k}(\xi)\nu^{6k},$$
$$\det \xi(\nu) = \sum_{\footnotesize{\begin{array}{l}k=-d\\k \, \hbox{is odd}\end{array}}}^{d}q_{k}(\xi)\nu^{3k},$$
where $p_{k}(\xi)$ and $q_{k}(\xi)$ are polynomials with respect to the elements of $\xi_{k}$
in $\xi(\nu)=\sum_{k=-d}^{d} \xi_{k}\nu^k$. Moreover, $p_{-k}=\bar{p_{k}},
 q_{-k}=-\bar{q_k}$. By Riemann-Hurwitz formula ${\stackrel{\circ}{g}}=6d-2$.

The quotient ${\hat \Gamma} = {\stackrel{\circ}{\Gamma}}/\gamma$ plays central role for explicit construction bellow. The three-sheeted covering
 $\stackrel{\circ}{\Gamma} \rw {\hat \Gamma}, 
 (\mu,\nu)\mapsto (\mu,\lambda),\nu^{3}=\lambda$ 
is unramified, where $\stackrel{\circ}{\Gamma}$ may be thought of as
the Riemann surface of the function $\nu = \sqrt[3]{\lambda}$ on $\hat \Gamma$.
The holomorphic involution $\sigma : (\mu, \lambda) \to (-\mu, -\lambda)$ and the
anti-holomorphic involution $\tau: (\mu, \lambda) \to (-\bar{\mu}, \f{1}{\bar{\lambda}})$
act on $\hat \Gamma$. By Riemann-Hurwitz formula ${\stackrel{\circ}{g}}= 3{\hat g}-2$ hence
${\hat g}=2d$.

On $\stackrel{\circ}{\Gamma}$ there exist three infinity points and three zero points, each 
situated on its copy of $\hat\Gamma$. We shall denote these points by $\infty^{I},\infty^{II},\infty^{III}$ and $0^{I},0^{II},0^{III}$ respectively. The automorphism $\gamma$ circularly acts on these infinity points
$\gamma(\infty^{I})=\infty^{III}, \gamma(\infty^{II})=\infty^{I},\gamma(\infty^{III})=\infty^{II}$ and $\gamma$ acts on these zero points similarly. Denote their projection on $\hat\Gamma$ by $P_{\infty}$ and $P_0$ respectively. Hence ${\hat \Gamma}\rw {\mathbb CP^1}$ is a three-sheeted cover of the $\lambda$-plane, and over $P_{\infty}$ and $P_0$ all sheets are glued, that is, the function 
$\lambda(P)$ on $\hat \Gamma$ having a pole of order 3 at $P_{\infty}$ and a zero of order 3 at $P_0$. We should note that although $\hat \Gamma$ is not hyperelliptic it shares many properties with one in the CMC tori. Such algebraic curve has its separating interests.  

Denote the points of $\hat \Gamma$ by $P=(\mu,\lambda)$. Since $\sigma$ leaves two distinct points $P_{\infty}$ and $P_0$ on $\hat \Gamma$ fixed, ${\hat \Gamma} \rw \Gamma={\hat \Gamma}/\sigma$ is a two-sheeted covering with two branched points. 
Therefore by Riemann-Hurwitz formula
 ${\hat g}=2g$ is twice of the genus $g$ of the factor surface $\Gamma$ and $g=d$ clearly.

Due to 
$$\xi(\lambda)_{z}=[\xi(\lambda), U(\lambda)],$$
$$\xi(\lambda)_{\bar z}=[\xi(\lambda), V(\lambda)],$$
the system
\begin{equation*}
\phi_{z}=\phi U, \phi_{\bar z}=\phi V, \phi\xi=\mu\phi
\end{equation*}
has a common (row) vector valued solution $\phi(P,z,{\bar z})$, which is called the {\it Baker Akhizer function}. 
\begin{remark}In the classical theory of integrable systems, the Baker-Akhizer function is defined as a column vector. For fitting into our situation, here we modify it to a row vector.
\end{remark}

In the finite-gap integration theory of the Tzit\'eica equation \eqref{db} usually a gauge equivalent function (see \cite{CS},\cite{Sharipov})
\begin{eqnarray*}
\psi = \phi \left(
\begin{array}{rcl}
\nu^{-2}& & \\
 &\nu^{-1}e^{\f{\omega}{2}}& \\
 & &e^{-\f{\omega}{2}} 
\end{array}
\right),
\end{eqnarray*}
is used.
The values of $\psi$ on different sheets of the covering ${\hat\Gamma}\rw {\mathbb CP^1}$
substituted in corresponding rows form the matrix $\Psi$ satisfying
\begin{eqnarray}\label{shar1}
{\Psi}_{z}={\Psi}\left(\begin{array}{ccc}
                                               0&0&{\bf i}\lambda\\
                                              {\bf i}&\omega_{z}&0\\
                                                0&{\bf i}&-\omega_z
                                                 \end{array} 
                                                    \right),
\end{eqnarray}
\begin{eqnarray}\label{shar2}
{\Psi}_{\bar z}={\Psi}\left(\begin{array}{ccc}
0&{\bf i}e^{\omega}&0\\
0&0&{\bf i}e^{-2\omega}\\
{\f{{\bf i}}{\lambda}}e^{\omega}&0&0
\end{array}
\right)
\end{eqnarray}
where ${\lambda}={\nu}^3$.  

After some computations one can prove the following analytic properties of $\psi$.
\begin{itemize}
\item[(1).]$\psi$ is a meromorphic function on ${\hat \Gamma} \backslash \{P_{0},P_{\infty} \}$. 
The pole divisor $\hat\mathcal{D}$ of $\psi$
on ${\hat\Gamma}\backslash \{P_{0},P_{\infty}\}$ is independent of $z,{\bar z}$ non-special divisor of degree $\hat g$. The Abel map of $\hat\mathcal D$ satisfies 
\begin{equation}\label{prymian}
{\mathcal A}({\hat\mathcal{D}}+\sigma {\hat\mathcal {D}})={\mathcal A}(P_{0}+P_{\infty}+ 
{\hat \mathcal C}),
\end{equation}
\begin{equation}\label{real}
{\mathcal A}({\hat\mathcal D}+\tau {\hat\mathcal D})={\mathcal A}(P_{0}+P_{\infty}+{\hat\mathcal  C}),
\end{equation}
where $\hat \mathcal C$ is the canonical class of $\hat \Gamma$.
\item[(2).] $\psi$ has essential singularities at the points $P_{0},P_{\infty}$ of the form
\begin{eqnarray*}
\psi_{1}(P)=\nu^{-3}e^{{\bf i}\nu z}(1+ o(1)),\\
\psi_{2}(P)=\nu^{-2}e^{{\bf i}\nu z}(1+ o(1)), \\
\psi_{3}(P)=\nu^{-1}e^{{\bf i}\nu z} (1+o(1))
\end{eqnarray*}
in the neighborhood of $P_{\infty}$, and 
\begin{eqnarray*}
\psi_{1}(P)=\nu^{-3}e^{{\f{{\bf i}}{\nu}}\bar z}(1+o(1)),\\
\psi_{2}(P)=\nu^{-2}e^{{\f{{\bf i}}{\nu}}\bar z}(e^{\omega}+o(1)),\\
\psi_{3}(P)=\nu^{-1}e^{{\f{{\bf i}}{\nu}}\bar z}(e^{-\omega}+o(1))\\
\end{eqnarray*}
in the neighborhood of $P_{0}$.
\end{itemize}

\section{Baker-Akhiezer function. Formulas}

The function $\psi$ is defined uniquely by its analytic properties and can be explicitly expressed in terms of 
Prym theta functions and Abelian integrals \cite{CS}. 
Let  $\hat \Gamma$ be a nonsingular Riemann surface of even genus ${\hat g}=2d, 
d\equiv 1\bmod 6$ with two distinguished point 
$P_{\infty}$ and $P_0$, on which there is a meromorphic function $\lambda(P)$ with divisor 
of zero and poles $3P_{0}-3P_{\infty}$ and on which two involutions are defined: 
the holomorphic involution $\sigma$, which acts in accordance with $\lambda(\sigma P)=-\lambda(P)$,
and an antiholomorphic involution $\tau$ of separating type such that 
$\lambda(\tau P){\overline {\lambda(P)}}=1$. 

Choose a canonical basis $\{a_{i},b_{i},i=1, \cdots, {\hat g}\}\in H_{1}({\hat \Gamma},{\mathbb Z})$ and
let $\omega_{1},\cdots , \omega_{\hat g}$ be the dual basis 
$$\oint_{a_i}\omega_{j} = 2\pi {\bf i} \delta_{ij}$$
of holomorphic differentials. The period matrix 
$$B_{ij}=\oint_{b_j}\omega_{i}$$
determines the Riemann theta function
$$\theta(z) = \sum_{N\in {\mathbb Z}^{\hat g}} \exp ({\f{1}{2}}\langle BN, N\rangle + \langle z, N \rangle), 
\quad z\in{\mathbb C^{\hat g}},$$
which is period with periods $2\pi {\bf i}{\mathbb Z^{\hat g}}$
\begin{equation}\label{period}
\theta(z+2\pi {\bf i} N)=\theta(z).
\end{equation}
Here the diamond brackets denote the Euclidean scalar product.
Let $\Lambda$ be the lattice generated by the vectors $2\pi {\bf i}e_{k}, Be_{k}, k=1,\cdots,{\hat g}$ 
where the vectors $e_{k}=(\delta_{k1},\cdots,\delta_{k{\hat g}})$. The complex torus
$J(\hat\Gamma)={\mathbb C^{\hat g}}/{\Lambda}$ is called the Jacobian of the Riemann surface $\hat \Gamma$. 
And the Abel mapping ${\mathcal A}(P) = ({\mathcal A}(P_1),\cdots,{\mathcal A}(P_{\hat g}))$,
$${\mathcal A}: {\hat\Gamma}\rw J(\hat\Gamma),$$
is defined by
$${\mathcal A}_{i}(P)=\int_{X}^{P}\omega_{i} \quad (i=1,\cdots, {\hat g}).$$
Here the basis point $X$ is often conveniently placed at the point $P_{\infty}$. 

For curves possessing involutions we can define the Prym variety  in the classical algebraic 
geometry (see, for example, \cite{Fay}). In our case, we can choose a basis 
$\{a_{i},b_{i},i=1, \cdots, {\hat g}\}$ such that 
$$\sigma a_{i}=-a_{i+g}, \quad \sigma b_{i} = -b_{i+g}. \quad i=1,\cdots,g.$$
Let $\varphi$ map $\mathbb C^g$ to $\mathbb C^{\hat g}$:
$$\varphi : (z_{1},\cdots, z_{g}) \rw (z_{1}, \cdots, z_{g},z_{1},\cdots, z_{g}).$$
With the help of the Abel mapping $\mathcal A$, the divisor map ${\hat\mathcal D}\mapsto \sigma {\hat\mathcal D}$ induces an involution 
$\sigma : J(\hat \Gamma) \rw J(\hat \Gamma)$. The subset 
$Prym(\hat\Gamma)=\{\Delta\in J(\hat\Gamma)|\sigma \Delta =-\Delta \}$ is called the {\it Prym variety} 
or the {\it Primian}
of the pair $({\hat\Gamma},\sigma)$. The Prymian is isomorphic to the torus ${\mathbb C^g}/{\Lambda}$,
where the lattice $\Lambda$ generated by the vectors $2\pi {\bf i} e_{k}$ and the column vectors of the matrix
$$\Pi_{kj} = \oint_{b_k} \omega_{j} + \omega_{j+g} \quad k,j = 1,\cdots,g.$$ 
The analogy form of the Abel mapping is 
\begin{equation}\label{prym}
{\mathcal U}_{j}(P) = \int_{X}^{P} (\omega_{j}+\omega_{j+g}), \quad j=1,\cdots, g.
\end{equation}
In addition, we can define the Prym theta function similarly as follows
$$\eta(z) = \theta(z|\Pi) = \sum_{N\in{\mathbb Z^g}}
\exp({\f{1}{2}} \langle \Pi N,N \rangle + \langle z,N \rangle ),
$$
where $z\in {\mathbb C^g}$. Suppose that ${\hat {\mathcal D}}$ is a fixed 
nonspecial divisor satisfying the above conditions \eqref{prymian}.
The Prym theta function constructed above solves the inversion problem with respect to \eqref{prym} 
in the class of divisors of degree $\hat g$ that satisfy condition \eqref{prymian} on $\hat \Gamma$ (see \cite{Fay}).
\begin{lemma}
If $e\in {\mathbb C^g}$ and $Q\in {\hat \Gamma}$, then either $F(P)=\eta({\mathcal U}(P-Q)-e) \equiv 0$
for all $P\in {\hat \Gamma}$ or the zeros of $F(P)$ form a divisor $\hat {\mathcal D}$ of degree $\hat g$ that satisfies
\eqref{prymian} and 
$${\mathcal A}({\hat {\mathcal D}}) = \varphi (e) +{\mathcal A}(Q-\sigma Q) +K,$$
where the constant vector $K$ is related to the canonical class $\hat \mathcal C$ of the surface $\hat \Gamma$ by 
the formula $2K={\mathcal A}({\hat\mathcal C} +P_{0} +P_{\infty})$. 
\end{lemma}

The lemma is proved in \cite{Fay} Corollary 5.6. By virtue of this lemma, we can specify 
the Baker-Akhizer function as follows \cite{CS}.

Choose local parameters $\nu^{-1}$ and $\varrho^{-1}$ in the neighborhood of the 
points $P_{\infty}$ and $P_{0}$, respectively, by means of the conditions
$$\nu^{3}=\lambda,\quad {\overline {\nu(\tau P)}} = \varrho(P),$$
such that $\nu^{-1}(P_{\infty})=0 = \varrho^{-1}(P_0)$ and $\nu^{-1}(\sigma P) = -\nu^{-1}(P)$, 
$\varrho^{-1}(\sigma P) = - \varrho ^{-1}(P)$. 

Let us introduce two abelian differentials of the second kind $\Omega_{\infty}, \Omega_{0}$ normalized by the conditions 
$$\oint_{a_i}\Omega_{\infty}=\oint_{a_i}\Omega_{0}=0,\quad i=1,\cdots,g,$$
and the following asymptotic behavior at the poles:
$$\Omega_{\infty}=d\nu +\cdots \quad \hbox{at} \, P_{\infty}, \quad \Omega_{0}=d\varrho+\cdots \quad \hbox{at} \, P_{0}.$$
We denote their periods upon $b$-cycles by
$$U_{i} = \oint_{b_i}\Omega_{\infty},\: V_{i}=\oint_{b_i}\Omega_{0},\, i=1,\cdots,g.$$
The involution $\sigma^*$ maps ${\Omega_{\infty}} \mapsto -\Omega_{\infty}, 
\Omega_{0}\mapsto -\Omega_{0}$.

The Baker-Akhizer function $\psi_1$ possessing the above properties 1 and 2 except \eqref{real}
can be described in terms with such data $\{{\hat \Gamma}, {\hat \mathcal D}\}$ uniquely up to a factor independent of $P$:
$$\psi_{1}(P)={\f{\eta({\mathcal U}(P)+{\bf i}Uz+{\bf i}V{\bar z}-e)}{\eta({\mathcal U}(P)-e)}}
\lambda^{-1}(P)\exp(\int_{P_\infty}^{P}{\bf i}z\Omega_{\infty}+{\bf i}{\bar z}\Omega_{0}).$$
The two remaining Baker-Akhiezer functions  $\psi_{2},\psi_{3}$ cannot be represented 
by such simple formulas in terms of Prym theta functions similarly. 
However, they can be computed by using the equations
$$\psi_{1z} = {\bf i}\psi_{2}, \quad \psi_{1\bar z} = {\f{\bf i}{\lambda}}e^{\omega}\psi_{3}.$$
But in fact we shall find in the following that it is unnecessary to specify the explicit formulas 
of $\psi_{2},\psi_{3}$. 
By the standard technique in the algebro-geometric theory of nonlinear integrable equation, the solution of the Tzitz\'eica equation \eqref{db} 
can be written as 
\begin{equation}\label{solu}
e^{\omega}=2\partial^{2}_{z\bar z}\ln \eta({\bf i}Uz+{\bf i}V{\bar z}-e)+c,
\end{equation}
where $c$ is the value of $\Omega_0$ at $\lambda=\infty$ satisfying $\Omega_0=cd\varrho+\cdots$ and $e$ satisfies ${\mathcal A}({\hat {\mathcal D}})=\varphi(e)+K$.

>From now on we put $z$ and $\bar z$ be mutually conjugated and consider the anti-holomorphic 
involution $\tau:\lambda\rw {\f{1}{\bar\lambda}}$. Then the Riemann surface $\hat \Gamma$ is divided by the fixed points of the involution $\tau$ (which are called the real ovals of $\tau$) into two regions: the region $\hat\Gamma_0$, which contains the point $P_0$, and the region $\hat\Gamma_{\infty}$, which contains the point $P_{\infty}$. The number of real ovals, $|\lambda|=1$, does not exceed three. It is determined by the number of real tori on the Jacobian $J({\hat\Gamma})$, each of which contains of classes of divisors $\hat\mathcal D$ of degree $\hat g$ satisfying the condition \eqref{real}. By virtue of \eqref{real}, every real divisor $\hat\mathcal D$ determines a certain Abelian differential 
$\alpha(P)$ of the third kind with zeros at the points of the divisor ${\hat\mathcal D}+\tau{\hat\mathcal D}$ and a pair of simple poles at $P_0$ and $P_{\infty}$ with residues $+i$ and $-i$, respectively. There exists one nonsingular torus $T_0$ which is distinguished among the other real tori by the fact that for the divisors from 
$T_0$ the differential $\alpha(P)$ is positive on all ovals of $\tau$ with respect to the natural orientation on the boundary of $\hat\Gamma_{\infty}$. Having fixed the torus $T_0$, we consider its subset consisting of the divisors satisfying the condition \eqref{real}, which is a real torus in the Prymian Prym$(\hat\Gamma)$ of the Riemann surface $\hat\Gamma$. We still denoted it by $T_0$.

Now let us consider the Riemann surface $\hat \Gamma$ with all branch points $\lambda_{i},
-\lambda_{i},i=1,\cdots,g$ divided into pairs ($|\lambda_{i}|\neq 1$)
$${\bar \lambda}_{2i-1}=\lambda_{2i}^{-1},\quad i=1,\cdots, g.$$
Canonical basis of cycles can be choose such that $\tau$ acts on it as follows 
$$\tau a_{i} = -a_{i}, \quad \tau b_{i} = b_{i}-a_{i}+\sum_{k=1}^{g}a_{k},$$
$$\tau a_{i+g} = -a_{i+g},\quad \tau b_{i+g} = b_{i+g}-a_{i+g}+\sum_{k=1}^{g}a_{k+g}.$$
Then the dual holomorphic differentials satisfy $\tau^{*}\omega_{i} = {\bar \omega}_{i}$ for $i=1,\cdots,g$
and the Prym period matrix $\Pi$ satisfies
\begin{eqnarray*}
  {\bar \Pi}_{kj}&=& \oint_{b_k}{\bar \omega}_{j} + {\bar\omega}_{j+g}
     =\oint_{b_k}\tau^{*}\omega_{j}+\tau^{*}\omega_{j+g} \\
     &=&\oint_{\tau b_k}\omega_{j}+\omega_{j+g}\\
     &=&\Pi_{kj}-2\pi {\bf i} \delta_{jk}+2\pi {\bf i},
\end{eqnarray*}
{\it i.e.,}
$${\bar\Pi}=\Pi+2\pi {\bf i} ({\bf 1}-I),$$
where ${\bf 1}_{kj}=1,I_{kj}=\delta_{kj}$. Therefore the Prym theta function satisfies the 
following conjugation condition
$${\overline {\eta(z)}}=\eta(\bar z),$$
which is easily checked using the evident identity
\begin{eqnarray*}
& &\exp\{{\f{1}{2}}\langle 2\pi {\bf i} ({\bf 1}-I)N,N \rangle\}\\
&=& \exp\{2\pi {\bf i} \sum_{k>j}N_{k}N_{j}\}=1.
\end{eqnarray*} 
Next we have 
$$\tau^{*}\Omega_{\infty}={\bar \Omega}_0,$$
hence the periods of $\Omega_{\infty}$ and $\Omega_0$ are mutually conjugated
$${\bar U}=V.$$
Therefore the vector ${\bf i}Uz+{\bf i}V{\bar z}-e$ is imaginary when ${\hat {\mathcal D}\in T_0}$.

\begin{theorem}
If ${\mathcal D}\in T_0$, then $\omega(z,\bar z)$ is a real nonsingular solution of \eqref{db}.  
\end{theorem}  
In relation to the functions $\psi_{1}, \psi_{2}, \psi_{3}$, the condition ${\hat \mathcal D}\in T_{0}$ gives 
\begin{eqnarray}\label{conjugate}
\psi_{1}(\sigma P) &=& -\lambda^{-2}(P){\overline{\psi_{1}(\tau P)}},\no\\
\psi_{2}(\sigma P) &=& -\lambda^{-1}(P)e^{\omega}{\overline{\psi_{3}(\tau P)}},\\
\psi_{3}(\sigma P) &=& \lambda^{-1}(P)e^{-\omega}{\overline{\psi_{2}(\tau P)}}.\no
\end{eqnarray} 
For the spectral problem \eqref{shar1} - \eqref{shar2}, we can define the pairing as follows 
which due to Sharipov :
\begin{eqnarray}\label{pair}
& & \Omega(P,Q) = \{ \psi(P)|\psi(\sigma Q)\}\\
&=& \{\psi_{3}(P)\psi_{2}(\sigma Q)\lambda (P) - \psi_{2}(P)\psi_{3}(\sigma Q)\lambda (P)
- \psi_{1}(P)\psi_{1}(\sigma Q)\lambda^{2}(P)\}.\no
\end{eqnarray}
Differentiating \eqref{pair} with respect to $z$ and $\bar z$, and taking into account \eqref{shar1} and 
\eqref{shar2}, we obtain the relations
$$\partial_{z}\Omega(P,Q) = {\bf i}[\lambda(Q) - \lambda(P)] \lambda(P) \psi_{2}(P)\psi_{1}(\sigma Q),$$
$$\partial_{\bar z}\Omega(P,Q)= {\bf i} e^{\omega}[\lambda(P)\lambda^{-1}(Q)-1]\lambda(P)\psi_{1}(P)\psi_{3}(\sigma Q),$$
from which we can see that when $P=Q$ the function \eqref{pair} does not depend on $z$ and $\bar z$. 
Moreover, the function $W(P) = \Omega(P,P)$ is meromorphic on $\hat \Gamma$ and can be calculated 
explicitly:
$$W(P) = {\f{{\bf i}d\lambda(P)}{\lambda(P)\alpha(P)}}.$$
Because the covering $\lambda : {\hat\Gamma}\rw {\mathbb C}$ is a three-sheeted one, each value of the function
$\lambda (P)$ is attained with multiplicity three, i.e., at three different points $P_{1},P_{2}$, and $P_{3}$. 
$\Omega(P,Q)$ has the following resonance property:
\begin{equation}\label{reson}
    \Omega(P_{i},P_{j})= 
        \left\{
          \begin{array}{cl}
              W(P_i) &  \mbox{for } P_{i}=P_{j},\\
              0   &  \mbox{for } P_{i}\neq P_{j},
\end{array}
\right.
\end{equation}
which in fact follows from that eigenvectors associated different eigenvalue of a normal matrix Hermitian orthogonal to each other. For every value of $\lambda \in S^{1}$, the points $P_{1}, P_{2}$, and $P_{3}$ lie on the ovals 
of the anti-holomorphic involution $\tau$ and are unchanged by the action of $\tau$. By them we form a matrix
\begin{eqnarray*}
\Psi = \left(\begin{array}{rcl}
\f{\psi_{1}(P_1)}{\sqrt{W(P_1)}}
&\f{e^{-\f{\omega}{2}}\psi_{2}(P_1)}{\sqrt{W(P_1)}}
&{\f{e^{\f{\omega}{2}}\psi_{3}(P_1)}{\sqrt{W(P_1)}}}\\
\f{\psi_{1}(P_2)}{\sqrt{W(P_2)}}
&{\f{e^{-\f{\omega}{2}}\psi_{2}(P_2)}{\sqrt{W(P_2)}}}
&{\f{e^{\f{\omega}{2}}\psi_{3}(P_2)}{\sqrt{W(P_2)}}}\\
{\f{\psi_{1}(P_3)}{\sqrt{W(P_3)}}}
&{\f{e^{-\f{\omega}{2}}\psi_{2}(P_3)}{\sqrt{W(P_3)}}}
&{\f{e^{\f{\omega}{2}}\psi_{3}(P_3)}{\sqrt{W(P_3)}}}
\end{array}
\right).
\end{eqnarray*}
>From \eqref{conjugate} and the invariance of the points $P_{1},P_{2}$ and $P_{3}$ with respect to $\tau$, 
 the resonance property \eqref{reson} leads to the relation
$$\psi_{1}(P_i){\overline{\psi_{1}(P_{j})}}+e^{-\omega}\psi_{2}(P_i){\overline{\psi_{2}(P_j)}}
+e^{\omega}\psi_{3}(P_i){\overline{\psi_{3}(P_j)}}=W(P_i)\delta_{ij},$$
which is equivalent to the matrix $\Psi$ is a unitary matrix. From this relation there also follow the 
reality and non-negativity of the values of the function $W(P)$ on the ovals of $\tau$, indicating 
the radicals in the matrix $\Psi$ are real. The columns of this unitary matrix $\Psi$ form a frame in $\mathbb C$
which solve the gauged spectral problem 
\begin{eqnarray*}
\Psi_{z}=\Psi \left(\begin{array}{rcl}
0&0&{\bf i}\lambda e^{\f{\omega}{2}}\\
{\bf i}e^{\f{\omega}{2}}& \f{\omega_z}{2}&0\\
0&{\bf i}e^{-\omega}&-\f{\omega_z}{2}
\end{array}
\right),\quad
\Psi_{\bar z}=\Psi\left(\begin{array}{rcl}
0&{\bf i}e^{\f{\omega}{2}}&0\\
0&-\f{\omega_{\bar z}}{2}&{\bf i}e^{-\omega}\\
{\f{\bf i}{\lambda}}e^{\f{\omega}{2}}&0&\f{\omega_{\bar z}}{2}
\end{array}
\right).
\end{eqnarray*}
Passing a gauged transformation $\sigma_{0}=\Psi \hbox{diag}(\nu^2, \nu, 1)$,
we have
\begin{theorem}
The $SU(3)$-valued solution $F(z,{\bar z}, \nu)$ of the linear system  \eqref{stan}, \eqref{1}, \eqref{2}
 with the coefficient \eqref{solu} 
is given by 
\begin{eqnarray*}
F&=&\left(\begin{array}{rcl}
\nu^{2}{\f{\psi_{1}(P_1)}{\sqrt{W(P_1)}}}
&\nu{\f{e^{-\f{\omega}{2}}\psi_{2}(P_1)}{\sqrt{W(P_1)}}}
&{\f{e^{\f{\omega}{2}}\psi_{3}(P_1)}{\sqrt{W(P_1)}}}\\
\nu^{2}{\f{\psi_{1}(P_2)}{\sqrt{W(P_2)}}}
&\nu{\f{e^{-\f{\omega}{2}}\psi_{2}(P_2)}{\sqrt{W(P_2)}}}
&{\f{e^{\f{\omega}{2}}\psi_{3}(P_2)}{\sqrt{W(P_2)}}}\\
\nu^{2}{\f{\psi_{1}(P_3)}{\sqrt{W(P_3)}}}
&\nu{\f{e^{-\f{\omega}{2}}\psi_{2}(P_3)}{\sqrt{W(P_3)}}}
&{\f{e^{\f{\omega}{2}}\psi_{3}(P_3)}{\sqrt{W(P_3)}}}
\end{array}
\right)\\
&:=&(y,F_{1},F_{2}).
\end{eqnarray*} 
Thus $x=[y(z, {\bar z}, 1)]:M \rw {\mathbb CP^2}$ gives a totally real superconformal minimal surface with the 
metric $g=2e^{\omega}dzd{\bar z}$ and $\psi=-1$.
\end{theorem}
As an evident corollary of \eqref{period}, the periodicity conditions can be given as follows:
\begin{theorem}
The immersion is doubly-periodic if the lattice is generated by the vectors $Z_{1},Z_{2}$
such that 
$${\mathrm Re} (Z_{k}U)\in \pi{\mathbb Z}^{g}, \quad 
{\mathrm Re} ({Z_k}\int_{P_{\infty}}^{P}\Omega_{\infty})\in\pi{\mathbb Z}, \quad k=1,2,$$
where $P=(\mu,\lambda)$ is chosen on the fixed ovals of $\hat \Gamma$.
\end{theorem}

It follows from standard techniques \cite{Bobenko91} that \eqref{solu} will yield the area formula 
of the totally real superconformal minimal torus. A detailed study of this aspect will be done elsewhere.


\begin{acknow}
The authors are grateful to Prof. Y. Ohnita
who initiated our research in this field during his visit to Peking university 
in the autumn of 1998 and revised carefully our preprint during the second international conference on harmonic morphisms and harmonic maps at the Centre International de Rencontres 
Mathematiques, Luminy in 2001. 
The authors are greatly indebted to Professors A.I. Bobenko, H.Z. Li,
A.V. Mikhailov,R.A. Sharipov, C.P. Wang and Y.J. Zhang for their assistance and useful comments,
and Dr. D. Joyce for his interest in this work.
The first author would like to express her sincere thanks to Prof. W.H. Chen for his continuous encouragements and helpful conversations.
Most of the results of this paper was reported by the first author during the second international conference on harmonic morphisms and harmonic maps at the CIRM in 2001. 
She would like to thank the organizers for their hospitality and the participants for their interests and helpful comments.

The first author is partially supported by the project No. 19871001 of NNSF of China. The second author is partially supported by the project No. G1998030601 of NKBRSF of China.
\end{acknow}
 

\end{document}